\documentclass[runningheads]{llncs}
\usepackage{multicol}

\newcommand\lgr{\longrightarrow}

\newcommand\cala{\mathcal{A}}
\newcommand\cals{\mathcal{S}}
\newcommand\calc{\mathcal{C}}
\newcommand\calt{\mathcal{T}}
\newcommand\calr{\mathcal{R}}

\begin{document}

\titlerunning{A Brownian quasi--helix 
 in $\bbbr^4$ built from an automatic sequence}
\authorrunning{J.-P. {Kahane}}

\title{A Brownian quasi--helix\\ in $\bbbr^4$, built from an automatic sequence}

\vskip6mm

\subtitle{S\'eminaire sur l'interface\\ entre l'analyse harmonique et la th\'eorie des nombres\\
\vskip2mm
Luminy 21 octobre 2005}
\vskip2mm

\author{J.-P. \textsc{Kahane}}

\institute{Universit\'e Paris--Sud, Laboratoire de Math\'ematique, B\^atiment 425,\\
F--91405 Orsay Cedex, France\\
}


\maketitle

\section{History and definitions}

The history of Brownian motion was described a number of times. In 1905, Einstein established his celebrated formula
$$
\overline{\Delta X^2} = \frac{RT}{N}\frac{1}{3\pi\mu a} \Delta t
$$
for spherical particles of radius $a$ suspended in a liquid of viscosity $\mu$ at temperature $T$ ; the first member, $\overline{\Delta X^2}$, is the average of the squares of their displacements during an interval of time $\Delta t$ ; $R$ is the constant of perfect gaz, and $N$ the Avogadro number. In the following years Jean~Perrin made a series of experiments leading to a new determination of the Avogadro number, and observed that the very irregular motion of particles resembled the nowhere differentiable  functions of mathematicians. Norbert Wiener introduced what he called ``the fundamental random function''  as a mathematical model for the physical Brownian motion. It was called immediately ``the Wiener process'', and later on, following Paul L\'evy, ``the Brownian motion''. Wiener gave several versions of the construction and derived a number of  fundamental properties, L\'evy developed the theory to a high point of sophistication, and it is now a mathematical object of common use as well as a mine of interesting problems.

Here is the theory as it appears from the last exposition made by Norbert Wiener. The problem is to construct a random process $X(t,\omega)$, also denoted by $X(t)$ $(=X(t,\cdot))$, ($t$ the time, $\omega\in\Omega$ the probability space) such that

\vskip2mm

1) for almost all $\omega$ $X(t,\omega)$ is a continuous function of $t$

2) $X(t)$ is a Gaussian process, meaning that the distribution of any $n$--uple $(X(t_1),X(t_2),\ldots X(t_n))$ is Gaussian

3) this Gaussian process has stationary increments, meaning that the distribution of $X(t)-X(s)$ depends on $t-s$ only

4) it satisfies a normalized Einstein equation, that is
$$
 ||X(t) -X(s) ||_2^2 = |t-s|
$$
where the norm is taken in $L^2(\Omega)$.


\noindent Here is such a construction. Let $\mathcal{H}$ be an infinite--dimensional subspace of $L^2(\Omega)$ consisting of Gaussian centered variables, and $W$ an isometric linear mapping of $L^2(I)$ ($I=\bbbr$, or $\bbbr^+$, or $[0,1]$) into $\mathcal{H}$. Let $\chi_t$ be the indicator function $1_{[0,t]}$. Then
$$
X(t) = W (\chi_t)
$$
satisfies all conditions 2) to 4). Moreover, given an orthonormal basis of $L^2(I)$, $(u_n)$, its image by $W$ is a normal sequence (sequence of independent Gaussian normalized random variables $(\xi_n)$), and expanding $\chi_t$ in the form
$$
\chi_t = \Sigma\ a_n(t) u_n \qquad (\hbox{in }L^2(I))
$$
results in an expansion of $X(t)$ as a random series of functions :
$$
\chi(t)= \Sigma\ a_n(t) \xi_n \qquad (\hbox{in }L^2(\Omega))
$$
or, more explicitly,
$$
X(t,\omega) = \Sigma\ a_n(t) \xi_n(\omega)\,.
$$
To prove condition 1), it is enough to establish that the series in the second member converges uniformly in $t$ for almost all $\omega$, and this is done rather easily when the $u_n$ are classical orthonormal bases.

By definition, an helix is a curve in a Hilbert space, parametrized by $\bbbr$, such that the distance between two points depends only on the distance of the parameters :
$$
||X(t) - X(s)||_2^2 = \psi(t-s)\,,
$$
and $\psi(\cdot)$ is called the helix function. A translation of the parameter results in a isometric motion of the curve onto itself. It is the abstract model for all Gaussian processes with independent increments.

When $\psi(t)=|t|$ we say that the curve is a Brownian helix. In contrast with the realizations of the Brownian motions (the functions $t\lgr X(t,\omega)$ when $\omega$ is fixed), the Brownian helix is a very regular curve. However some basic properties of the Brownian motion can be read on the Brownian helix : its Hausdorff dimension is $2$, its $2$--dimensional Hausdorff measure is nothing but $dt$, and any three points on the curve are the vertices of a rightangle triangle : the increments starting from a point are orthogonal to the past (therefore, independent from the past).

Simple examples of helices are :

 1) the line $(\psi(t)=a^2t^2)$
 
 2) the circle $(\psi(t)=r^2\ \sin^2 \omega t)$
 
 3) the three--dimensional helices $(\psi(t) =a^2t^2 + r^2\ \sin^2 \omega t)$
 
 4) generalizations of those, with
$$
\psi(t) = a^2t^2 + \int\ \sin^2 \omega t\ \mu(d\omega)\,,
$$
\noindent
where $\mu$ is a positive measure on $\bbbr^+$ such that the integral is finite. Actually this is the general form of an helix function.

Except when $\mu$ is carried on a finite set, the helix cannot be imbedded in a finite dimensional Euclidean space.

At the end of the 70s, Patrice Assouad developed a theory of Lipschitz embeddings of a metric space into another \cite{assouad}. He introduced and built quasi--helices in Euclidean spaces, meaning that
$$
0<a<\frac{||X(t)-X(s)||_2^2}{\psi(t-s)} < b <\infty
$$
for some $a$ and $b$ and all $t$ and $s$. When $\psi(t)=|t|$ we call them Brownian quasi--helices.  Assouad constructed Brownian quasi--helices in Euclidean $\bbbr^n$ for $n\geq 3$, and this gives a new way to prove that the realizations of Brownian motion are continuous $a.s.$. He asked me the question whether $a$ and $b$ can be taken near $1$ when $n$ is large, that is, whether the Brownian helix can be approximated (in this sense) by Brownian quasi--helices. I gave a positive answer with an explicit construction, and it was published in my paper on Helices and quasi--helices \cite{kahane}.

\section{A construction of Brownian quasi--helices by means of Walsh matrices}

Let us consider $\bbbr^{2^n}$ $(n\geq 1)$ as a Euclidean space. Let $N=2^n$. If we want to construct a function $X:\bbbn \lgr \bbbr^N$ such that $X(0)=0$ and $||X(t) -X(s)||^2 = |t-s|$ when $|t-s| \leq N$ we have to choose an orthonormal basis $u_0,u_1,\ldots u_{N-1}$, define $u_{N+j}=u_j$, and write
$$
X(t) = \sum_{0\leq j\leq t-1} \pm u_j\,.
$$
At this stage there is no restriction on the signs $\pm$, and we may choose $+$ when $0\leq j\leq N-1$. If we try to obtain $||X(2t)-X(2s)||^2 = 2|t-s|$, $||X(4t) - X(4s)||^2 = 4|t-s|$ etc when $|t-s|\leq N$, whe have more and more conditions on the $\pm$ and we are led to the following construction.

We define the Walsh matrix of order $N$ as the $N\times N$ matrix obtained as the $n^{\rm th}$ tensor power of the matrix $\left(
\begin{array}{cc}
  1&   1   \\
  1&  -1   \\  
\end{array}
\right)$
, that is
 $$
\left(
\begin{array}{cc}
  1&   1   \\
  1&  -1   \\  
\end{array}
\right) \otimes \left(
\begin{array}{cc}
  1&   1   \\
  1&  -1   \\  
\end{array}
\right) \otimes \cdots \otimes \left(
\begin{array}{cc}
  1&   1   \\
  1&  -1   \\  
\end{array}
\right) \qquad (n\hbox{ times)}\,.
$$


\noindent
For example, the Walsh matrix of order 4 is
$$
M=M_2 = 
\left(
\begin{array}{cccc}
  1 &1   &1 &1   \\
  1& -1  & 1 &-1  \\
  1&+1   &-1 &-1 \\
  1 &-1 &-1 &1   
\end{array}
\right)
$$
and the matrix $M_{n+1}$ of order $2^{n+1}$ is obtained from $M_n$ as
$$
M_{n+1} = 
\left(
\begin{array}{cc}
  M_n&   M_n  \\
  M_n&  -M_n   \\  
\end{array}
\right)
$$
The $N^2$ first signs $\pm$ are those of the entries of the Walsh matrix, read line by line. In order to obtain the following signs, we extend the Walsh matrix by a series of vertical translations and change of signs of some lines according to the following rule : the first row is nothing but the whole sequence of entries, written from line to line and from left to right.

With this procedure we define $X(t)$ when $t$ is an integer and we can extend the construction to all $t>0$, then to all real $t$. It is proved in \cite{kahane} that we obtain a quasi--helix with $a$ and $b$ close to $1$ when $n$ is large enough : it is the answer to the question of Assouad.

However, it was not proved that the construction provides a quasi--helix when $n=2$ (it was remarked that it gives a Peano curve in the plane when $n=1$). The aim of the present paper is to give a detailed exposition of the case $n=2$ (most of it could be copied for $n>2$) and to prove that we obtain a quasi--helix. Instead of $t\in \bbbr$ we shall consider only $t\in \bbbr^+$ and a curve starting from $0$ $(X(0)=0)$. We shall investigate the geometric properties of the curve, some of them leading to open questions of a combinatorial or arithmetical nature.

The sequences that we construct are automatic in the sense of \cite{allou-shal}

\section{Description of the sequence}

\textbf{3.1} \hskip2mm It is a sequence of $+1$ and $-1$ as described before, in case $N=4$. We write it as a succession of $+$ and $-$ :
$$
++++\quad +-+- \quad ++-- \quad +--+ \quad ++++\ \cdots
$$
The gaps between the blocks of four letters have no meaning, except a help to understand the construction. The construction proceeds as follows : given the initial word of length $4^j$, we divide it into four words $A$, $B$, $C$, $D$ of equal length $4^{j-1}$ and write it $ABCD$ ; then
$$
A\ B\ C\ D\ A\ (-B)\ C\ (-D)\ A\ B\ (-C)(-D)\ A(-B)(-C)\ D
$$
is the initial word of length $4^{j+1}$. We shall give several equivalent definitions,  using substitutions, explicit expressions, or generating functions.

Beforehand let us write the sequence in a tabular form as in the previous section :
$$
\begin{tabular}{|cccc|rl}
\cline{1-4}
+ &+&+&+ & &\hskip 3cm $a_0$\ $a_1 $\ $a_2$ \ $a_3$\\
+ &$-$&+&$-$& \hskip 3mm $A$  & \hskip 3cm$a_4$\ $a_5$\ \hbox to 2.5cm{...........}           \\
+&+&$-$&$-$   &&\hskip 3cm\hbox to 4cm{....................}\\
+&$-$&$-$&+  & &\hskip 3cm.\hbox to 1.5cm{............}\ $a_{15}$\\
\cline{1-4}
+&+&+&+ &&\hskip3cm $a_{16}$ \hbox to 4cm{.............}\\
$-$&+&$-$&+ & \hskip 3mm $B$ &\hskip3cm \hbox to 4cm{....................}\\
+&+&$-$&$-$ &  &\hskip3cm \hbox to 4cm{....................}\\
$-$&+&+&$-$ & &\hskip 3cm.\hbox to 1.5cm{............}\ $a_{31}$\\
\cline{1-4}
+&+&+&+    &&\hskip 3cm $a_{32}$ \hbox to 4cm{.............}\\
+&$-$&+&$-$& \hskip 3mm $C$ &\hskip3cm \hbox to 4cm{....................}\\
$-$ &$-$& +&+ & &\hskip3cm \hbox to 4cm{....................}\\
$- $&+ &+&- & &\hskip3cm \hbox to 1.5cm{............}\ $a_{47}$\\
\cline{1-4}
+&+&+&+           &    &\hskip 3cm$a_{48}$ \hbox to 4cm{.............}\\
$-$ &+ &$-$ &+ &\hskip3mm$D$ &\hskip3cm \hbox to 4cm{....................}\\
$-$ &$-$ &+&+ &&\hskip3cm \hbox to 4cm{....................}\\
+ &$-$ &$-$ &+ &&\hskip3cm \hbox to 1.5cm{............}\ $a_{63}$\\
\cline{1-4}

&&&     &&\hskip 3cm$a_{64}$ \hbox to 4cm{.............}\\
&&&      &\hskip3mm$A$ &\hskip3cm \hbox to 4cm{....................}\\
&&&     &&\hskip3cm \hbox to 4cm{....................}\\
&&&     &&\hskip3cm \hbox to 1.5cm{............}\ $a_{79}$\\

\cline{1-4}

&&&     &&\hskip 3cm$a_{80}$ \hbox to 4cm{.............}\\
&&&      &\hskip3mm$-B$ &\hskip3cm \hbox to 4cm{....................}\\
&&&     &&\hskip3cm \hbox to 4cm{....................}\\
&&&     &&\hskip3cm \hbox to 4cm{....................}\\
\cline{1-4}
\end{tabular} 
$$
\vskip6mm
\noindent\textbf{3.2} \hskip2mm Let me give an explicit expression for $a_n$. Writing
$$
n= n_0 +4n_1 +\cdots +4^\nu n_\nu \qquad (n_\nu=1,2,3 ;\  n_j=0,1,2,3\ \hbox{if } j<\nu)
$$
the construction shows that
$$
a_n = a_{n_0+4n_1} \ a_m\,,\quad m=n_1 +4n_2+\cdots +4^{\nu-1}n_\nu\,,
$$
that is
$$
a_n = a_{n_0+4n_1} \ a_{n_1+4n_2} \cdots a_{n_{\nu-1}+4n_\nu}\,.
$$
In the second member we find $a_j$ $s$ with $j\leq 15$. Their value is $-1$ when $j=5,7,10,11,13$ and $14$ and $+1$ otherwise. Now let us express $a_n$ as a function of $n$ written  in the $4$--adic system of numeration. We obtain the formula
$$
a_n=(-1)^{A_n}
$$
\noindent
where $A_n$ is the number of $11,13,22,23,31,32$ in the $4$--adic expansion of $n$. For example, if $n= 1\ 3\ 2\ 0\ 0\ 1\ 1\ 1\ 0\ 2\ 3\ 1\ 1\ 1\ 2\ 2$ the significant links are
$$
1_- 3_-2\ 0\ 0\ 1_-1_-1\ 0\ 2_-3_-1_-1_-1\ 2_-2
$$
$A_n$ is nine and $a_n=-1$.

\vskip4mm

\noindent \textbf{3.3}  \hskip2mm Let me describe the sequence by means of a substitution rule.

We start from an alphabet made of eight letters : $+a,+b,+c,+d,-a,-b,-c,-d$. The substitution rule is
$$
\begin{array}{lll}
  +a& \lgr  &+a+b+c+d   \\
  +b&   &+a-b+c-d   \\
  +c&   &  +a+b-c-d\\
  +d &&+a-b-c+d\\
  -a && -a-b-c-d\\
  -b && -a+b-c+d\\
  -c &&-a-b+c+d\\
  -d &&-a+b+c-d\\ 
\end{array}
\leqno{(S_0)}
$$
The infinite word beginning with $+a$ and invariant under the substitution is
$$
W= +a+b+c+d+a-b+c-d +a+b-c-d+a-b-c+d+\cdots
$$
Replacing $a,b,c,d$ by $1$ (or, in a graphic way, in suppressing them), we obtain our sequence of $\pm1$ (or $\pm$).

\vskip4mm

\noindent \textbf{3.4} \hskip2mm Actually there is a simpler substitution rule leading to the same result, namely
$$
\begin{array}{lll}
  +a& \lgr  &+a+b   \\
  +b&   &+c+d   \\
  +c&   &  +a-b \\
  +d &&+c-d\\
  -a && -a-b\\
  -b && -c-d\\
  -c &&-a+b\\
  -d &&-c+d\\ 
\end{array}
\leqno{(S_1)}
$$
It can be checked immediately that $(S_1)(S_1)=(S_0)$.

\vskip4mm

\noindent \textbf{3.5} \hskip2mm The generating function of the sequence $(a_n)$ is 
$$
f(z) =a_0 +a_1 z +a_2 z^2+\cdots
$$
It can be defined using partial sums of order $4^n$. Let us introduce the matrix
$$
M(z) = 
\left(
\begin{array}{cccc}
  1&z   &z^2  &z^3 \\
 1 &-z   &z^2   &-z^3 \\
  1& z  &-z^2  &-z^3  \\
  1& -z  &-z^2  &z^3
\end{array}
\right)
$$
and define four sequences of polynomials by the formulas
$$
\left(
\begin{array}{c}
 P_0    \\
  Q_0   \\
R_0  \\
T_0   
\end{array}
\right) = 
\left(
\begin{array}{c}
 1    \\
 1   \\
1 \\
1   
\end{array}
\right) \qquad 
\left(
\begin{array}{c}
 P_{n+1}    \\
  Q_{n+1}   \\
R_{n+1}  \\
T_{n+1}   
\end{array}
\right) = M(z^{4^n}) 
\left(
\begin{array}{c}
 P_n    \\
  Q_n   \\
R_n  \\
T_n  
\end{array}
\right) 
$$
$(n=0,1,\ldots)$. Then
$$
\left(
\begin{array}{c}
 P_n    \\
  Q_n   \\
R_n  \\
T_n   
\end{array}
\right) =M(z^{4^{n-1}}) M(z^{4^{n-2}})\cdots M(z^4)M(z)
\left(
\begin{array}{c}
 1    \\
1   \\
1  \\
1   
\end{array}
\right) \,.
$$
When $|z|=1$, we have $M(z)M(\overline z)=4I$, therefore the matrix ${1\over2}M(z)$ is unitary, hence
$$
\begin{array}{ll}
|P_n|^2 + |Q_n|^2 + |R_n|^2 + |T_n|^2 &=4(|P_{n-1}|^2 + |Q_{n-1}|^2 + |R_{n-1}|^2 + |T_{n-1}|^2)\\
   &=4^n (1+1+1+1)=4^{n+1}
\end{array}
$$
We obtain the generating function as
$$
f(z) = \lim_{n\rightarrow \infty} P_n(z)
$$

We can write the generating function in a more interesting form :
$$
f(z) =f_0(z^4) +z f_1(z^4) +z^2 f_2(z^4) +z^3 f_3(z^4)\,,
$$
where the coefficients of the power series $f_0,f_1,f_2,f_3$ are the columns of the table in 3.0. In order to obtain these coefficients, we can start from $W$ in 3.2 and replace $a,b,c,d$ by $1,1,1,1$ (for $f_0)$, $1,-1,1,-1$ (for $f_1$), $1,1,-1,-1$ (for $f_2$) and $1,-1,-1,1$ (for $f_3$). Then $f_0=f$. Writing
$$
F(z) =
\left(
\begin{array}{c}
  f_0 (z)   \\
 f_1 (z) \\
  f_2 (z) \\
   f_3 (z) 
\end{array}
\right)
$$
the functional equation of the generating functions of the columns is 
$$
F(z) = M(z) F(z^4)\,.
$$

\section{Description of the curve}

\textbf{4.1} \hskip2mm Let $u_0,u_1,u_2,u_3$ be an orthonormal basis of the Euclidean space $\bbbr^4$, and define $u_{j+4}=u_j$ $(j=0,1,\ldots)$. The partial sums of the series
$$
a_0u_0 +a_1 u_1 +a_2u_2 +\cdots
$$
(that can be obtained from $W$ in 3.2 by replacing $a,b,c,d$ by $u_0,u_1,u_2,u_3$) will be denoted by $S(n)$. Then
$$
S(n) = a_0u_0 +a_1u_1+\cdots +a_{n-1}u_{n-1} \in \bbbz^4\,.
$$
It is easy to check on the table in 3.0 that
$$
S(16n) =4S(n)\qquad (n\in \bbbn)\,.
$$
Moreover it is not difficult to see (we shall be more specific later) that
$$
||S(n)-S(m)||^2 \leq b\ |n-m|
$$
for some $b<\infty$ and all $n$ and $m$. This allows, first, to define $S(t)$ when $t$ is a binary number via a formula $S(16^\nu t) = 4^\nu S(t)$, then to check that
$$
\begin{array}{c}
 S(16 t) =4S(t)   \\
||S(t)-S(s)||^2 \leq b |t-s|   \\
\end{array}
$$
for such numbers, then to extend $S(\cdot)$ by continuity on $\bbbr^+$, and check the above formulas for all $t\geq 0$ and $s\geq 0$.

The curve we consider is the image of $\bbbr^+$ by $S(\cdot)$.

Clearly (changing $t$ into $16t$) the curve is invariant under an homothety of center $0$ and ratio $4$. Our main purpose is to prove that it is a Brownian quasi--helix. We shall point out some geometric properties first.

\vskip4mm

\noindent\textbf{4.2} \hskip2mm The matrix $M$ transforms $u_0,u_1,u_2,u_3$ into $u_0',u_1',u_2',u_3'$ :
$$
(u_0',u_1',u_2',u_3') = M (u_0,u_1,u_2,u_3)
$$
and the partial sums of order $n$ of the series $\Sigma a_ju_j'$ are the partial sums of order $4n$ of the series $\Sigma a_ju_j$. Therefore the equation
$$
S(4t)=M\ S(t)
$$
holds true when $t=n\in \bbbn$ and by extension for all $t\in\bbbr^+$.

It is easy to check that the eigenvalues of $M$ are $2$ and $-2$, and that
$$
M
\left(
\begin{array}{cc}
  1&   1   \\
  1&   0   \\
  1&   0   \\
  -1 &1   
\end{array}
\right) =
\left(
\begin{array}{cc}
  2&   2   \\
  2&   0   \\
  2&   0   \\
  -2 & 2 
\end{array}
\right)\,, \qquad
M \left(
\begin{array}{cc}
  1&  0   \\
  -1&  1   \\
  -1&   -1   \\
  -1 & 0   
\end{array}
\right) = 
\left(
\begin{array}{cc}
  -2&   0   \\
  2&   -2   \\
  2&   2   \\
  2 & 0 
\end{array}
\right)\,.
$$
Let
$$
\begin{array}{cc}
  v_0 = {1\over2}(u_0+u_1+u_2-u_3)\,,&  v_1 = {\sqrt{2}\over2}(u_0+u_3)\hfill   \\
   v_2 = {1\over2}(u_0-u_1-u_2-u_3)\,, &    v_3 = {\sqrt{2}\over2}(u_1-u_2)\,.    \\
\end{array}
$$
They constitute an orthonormal system. The vectors $v_0$ and $v_1$ generate a plane, $P$, which is the eigenspace of the eigenvalue $2$, and $v_2$ and $v_3$ a plane, $Q$, corresponding to the eigenvalue $-2$. Expressed via the orthonormal basis $v_0,v_1,v_2,v_3$, the operator $M$ takes the form
$$
M'=U\ M \ U^{-1}= 2 
\left(
\begin{array}{cccc}
  1\     &0   &0    &0   \\
  0\     &1   &0    &0   \\
  0\     &0   &-1    &0   \\
  0\     &0   &0    &-1 
\end{array}
\right)\,,
$$
where $U$ is the unitary matrix carrying $u_0,u_1,u_2,u_3$ onto $v_0,v_1,v_2,v_3$. It means that the transformation $S(t) \lgr S(4t)$ is the product of a homothety of centre $0$ and ratio $2$ and an orthogonal symmetry with respect to the plane $Q$.

\vskip4mm

\noindent \textbf{4.3} \hskip2mm Clearly $M^2=4I$ ($I$ being the identity matrix), In turn, $M$ is the square of another simple matrix, as it can be guessed  from 3.2 and 3.3. Let us define
$$
T = \left(
\begin{array}{cccc}
  1\ &   0&    1    &0   \\
  1\ &   0&    -1    &0   \\
   0\ &   1&    0   &1 \\
   0 \ &  1&   0  &-1  
\end{array}
\right)\,.
$$
Then $M=T^2$.

The eigenvalues of $T$ are $\sqrt{2}$, $-\sqrt{2}$, $i\sqrt{2}$ and $-i\sqrt{2}$. The vectors
$$
w_0 = \frac{\sqrt{2}}{2} (v_0+v_1)\,, \qquad w_1 = \frac{\sqrt{2}}{2} (v_0-v_1)
$$
are eigenvectors corresponding to $\sqrt{2}$ and $-\sqrt{2}$. Defining $w_2$ and $w_3$ in such a way that $w_0,w_1,w_2,w_3$ is a direct orthonormal basis, and $W$ being the unitary matrix carrying $u_0,u_1,u_2,u_3$ onto $w_0,w_1,w_2,w_3$, we can write
$$
WTW^{-1}= T' =\sqrt{2} 
\left(
\begin{array}{cccc}
  1& 0  &0   &0   \\
 0 &  -1 &0  &  0  \\
  0&   0& 0  &-1 \\
  0 &0 &1 &0
\end{array}
\right)\,.
$$
It means that $T'$ is decomposed into :

1) an homothety of centre $0$ and ratio $\sqrt{2}$

2) a rotation of $\frac{\pi}{2}$ of the orthogonal projection on $Q$

3) a symmetry with respect to $w_1$ of the orthogonal projection on $P$.

\vskip1mm
In the same way as we obtained the equation $S(4t)=MS(t)$, we now have
$$
S(2t) = T\,S(t)
$$
and we have just given the interpretation of the transformation $S(t) \lgr S(2t)$ as a product of simple transformations.

\vskip4mm

\noindent \textbf{4.4} \hskip2mm
We have investigated the properties of the transformations $S(t)\lgr S(16t)$, $S(t)\lgr S(4t)$, $S(t)\lgr S(2t)$ as products of homotheties and isometries. Now we shall look at the effect of a translation of $t$ by an integer. We are interested in differences $S(t)-S(s)$.

Let us begin with integers $m<n<16^k$. Let us divide the series $a_0u_0+a_1u_1+\cdots$ into consecutive blocks of length $16^k$, so that the series reads
$$
+A + B +C +D +A -B +C -D + \cdots
$$
If $j=j_0 +4j_1$ $(j_0 =0,1,2,3,j_1\in \bbbn)$, the $j$--th term is of type $A,B,C,D$ according to the value of $j_0$ and its sign is $a_j$. Therefore
$$
S(n+j\cdot 16^k) - S(m+j\cdot 16^k) = a_j (S(n+j_0 16^k) -S(m+j_0 16^k))\,.
$$

If $0<s<t<1$, we can approximate $s$ and $t$ by $m\, 16^{-k}$ and $n\, 16^{-k}$ and we obtain
$$
S(t+j) -S(s+j) = a_j(S(t+j_0) -S(s+j_0))
$$
$(j_0=1,2,3,4,\ j_0=j$ modulo $4$).

This expresses that all arcs $\cala_j = S([j,j+1])$ are isometric (actually translates or symmetric according to the value of $a_j$) of one of the arcs $\cala_0,\cala_1,\cala_2,\cala_3$ (according to the value of $j_0$). Using 4.2, this  holds true when we replace the $\cala_j$s by $\cala_j^\nu = S([j2^\nu,(j+1)2^\nu])$ whatever $\nu\in \bbbz$.

\section{It is a Brownian helix}

\vskip2mm 

\textbf{5.1} \hskip2mm What we have to prove is that, writing
$$
a = \inf_{0<s<t} \frac{||S(t)-S(s)||}{\sqrt{|t-s|}} \leq \sup_{0<s<t} \frac{||S(t)-S(s)||}{\sqrt{|t-s|}} = b\,,
$$
we have
$$
a>0\,, \qquad b < \infty\,.
$$
We can write as well
$$
a = \inf_{m<n}\frac{||S(n)-S(m)||}{\sqrt{|n-m|}}\,, \quad b=\sup_{m<n} \frac{||S(n)-S(m)||}{\sqrt{|n-m|}}
$$

\vskip2mm

\noindent \textbf{5.2} \hskip2mm The easy part is $b<\infty$.

Let us first assume $[m,n] = [j 2^k, (j+1) 2^k]$. Then, according to 4.3,
$$
||S(n) - S(m)|| = 2^{k/2}\,.
$$
In the general case, let us decompose $[m,n]$ into such intervals in a minimal way, so that there are at most two intervals of the same length in the decomposition. If the largest length is $2^k$, we obtain
$$
||S(n)-S(m)|| \leq 2(2^{k/2} +2^{(k-1)/2} +\cdots ) \leq 2(1+\sqrt{2}) 2^{k/2}
$$
therefore
$$
||S(n) - S(m)|| \leq 2(1+\sqrt{2}) |n-m|^{1/2}\,.
$$
This gives
$$
b\leq 2(1+\sqrt{2})
$$

\noindent\textbf{5.3} \hskip2mm To prove $a>0$ is more tricky.

We shall use two lemmas.

\vskip2mm

\begin{lemma}  {There exists $\alpha>0$ such that} 
$$
||S(n+h) - S(n)|| \leq 1-\alpha
$$
{for all $n\in\bbbn$ and $h\in [-\frac{1}{2},\, \frac{1}{2}]$ (with $n+h\geq 0$).}
\end{lemma}

\vskip2mm

\begin{lemma} There exists an integer $A$ such that
$$
||S(n) -S(m)|| \geq 2
$$
for all integers $m$ and $n$ such that $n-m\geq A$.
\end{lemma}
\vskip2mm

Assuming that this is correct, the result is at hand : given $t$ and $s$ such that $t-s\geq A+1$, we can write $s=m+h$ and $t=n+h'$ with $h$, $h'\in [-\frac{1}{2},\frac{1}{2}]$ and $m-n\geq A$, therefore
$$
||S(t)-S(s)||\geq 2\alpha\,.
$$
Whenever $(A+1)2^k < t-s \leq (A+1)2^{k+1}$ $(k\in\bbbz)$, we have
$$
||S(t)-S(s)|| \geq 2\alpha 2^{k/2} \geq \frac{2\alpha}{\sqrt{2(A+1)}}\ |t-s|^{1/2}
$$
therefore
$$
a\geq \frac{2\alpha}{\sqrt{2(A+1)}}\,.
$$

\vskip4mm
\noindent\textbf{5.4} \hskip2mm Proof of Lemma 1.

From now on it may be useful to represent $S(n)$ on the table of 3.0, and also the differences $S(n)-S(m)$, as figures consisting of consecutive lines plus or minus part of a line above and below, in such a way that each column in the figure has a sum equal to the corresponding coordinate of $S(n)$ or $S(n)-S(m)$.

$$
\begin{array}{c}
\begin{tabular}{|cccc|}
\hline
+ &+&+&+ \\
+ &$-$&+&$-$       \\
+&+&$-$&$-$   \\
+&$-$&$-$&+\\
\cline{2-4}
\multicolumn{1}{|c|}{+} &&\\
\cline{1-1}
\end{tabular}
\\
S(17)
\\
\\
\\
\\
\\
\end{array}
\qquad
\begin{array}{c}
\begin{tabular}{|cccc|}
\hline
&&& \\
&&&       \\
&&&   \\
&&&   \\
\hline
&&& \\
&&&       \\
\cline{2-4}
\multicolumn{1}{|c|}{} &+&$-$ &$-$\\
\cline{1-1}
$-$ &+ &+ &$-$\\
\hline
\end{tabular}\\
\\
S(32) -S(25)
\end{array}
$$

Let us consider
$$
||S(16n+m)-S(16n)||^2
$$
for $n\in \bbbn$ and $m=0,\pm1,\pm2,\pm3,\ldots,\pm 8$. It is sufficient to consider the four cases $n=0,1,2,3$, and to look at the figures (depending on $m$) in each case. The result is
$$
\begin{array}{cc}
||S(16n+m)-S(16n)||^2 \leq 8  &(n\ \hbox{odd)}      \\
 ||S(16n+m)-S(16n)||^2 \leq 9 &(n\ \hbox{even)}       \\   
\end{array}
$$
with equality only when $m$ is odd (as for $S(32)-S(25)$). Therefore, going one step further,
$$
\begin{array}{ll}
||S\big(16n+m+\frac{P}{16}\big) - S(16n)|| &\leq \sqrt{8} +\frac{1}{4} \sqrt{9} \ \ (n\ \hbox{odd})   \\
\\
||S\big(16n+m+\frac{P}{16}\big) - S(16n)||  & \leq \sqrt{9} +\frac{1}{4} \sqrt{8} \ \ (n\ \hbox{even}) \\
\end{array}
$$
when $p=0,\pm1,\pm2,\ldots,\pm8$. Proceeding that way we finally obtain
$$
||S(16(n+t)) - S(16n)|| \leq 4(1-\alpha)\qquad \big(-{\textstyle{1\over2}} \leq t \leq{ \textstyle{1\over2}}\big)
$$
where
$$
4(1-\alpha)=\big(\sqrt{9}+{1\over4}\sqrt{8}\big)\big(1+{1\over16}+{1\over16^2}+\cdots \big) 
$$
as the second member is less than $4$ and that proves Lemma 1.

\vskip4mm

\noindent \textbf{5.5}  \hskip2mm  Proof of Lemma 2

Here again we look at the table. We can compute
$||S(n)-S(m)||$
when $n-m$ is in a given interval, and moreover give the couples $(m,n)$ for each the infimum is attained, and the expression of $S(n)-S(m)$ (that is, the coordinates with respect to $u_0,u_1,u_2,u_3$).

1) $4<n-m \leq 16$. It suffices to consider the first 16 lines of the table (fig.~1), since adding to $m$ and $n$ a multiple of 64 does not change $||S(n)-S(m)||$. We obtain
$$
\inf ||S(n)-S(m)||^2 =2
$$
realized for $(5,11)$, $(23,29)$, $(35,41)$ and $(53,59)$ :
\vskip2mm
$$
\begin{array}{cc}
 S(11) - S(5) = u_1-u_4\,,\hfill & S(29) -S(23) =u_2 -u_3\hfill  \\
  S(41) -S(35) = -u_2+u_3\,, &S(59) - S(53)=-u_1 +u_4   \\  
\end{array}
$$

\begin{center}
$\begin{array}{ll}
\begin{tabular}{|cccc|rr}
\cline{1-4}
&&&\\
\cline{2-4}
\multicolumn{1}{|c|}{} &$-$ &+&$-$ &$S(11)-S(5)$\ \  \\
\cline{1-1}\cline{4-4}
+&+&$-$&\multicolumn{1}{|c|}{}\\ 
\cline{1-3}
&&&\\
\cline{1-4}
&&&\\
\cline{4-4}
&&&\multicolumn{1}{|c|}{+} &$S(29)-S(23)$\\
\cline{1-3}
+ &+&$-$ &$-$\\
\cline{2-4}
\multicolumn{1}{|c|}{$-$} &&&\\
\cline{1-4}
&&&\multicolumn{1}{|c|}{+} &$S(41)-S(35)$\\
\cline{1-3} 
+ &$-$&+&$-$ \\
\cline{2-4}
\multicolumn{1}{|c|}{$-$} &&&\\
\cline{1-1}
&&&\\
\cline{1-4}
&&&\\
\cline{2-4}
\multicolumn{1}{|c|}{} &+&$-$&+ &$S(59)-S(53)$\\
\cline{1-1}\cline{4-4}
$-$&$-$ &+ &\multicolumn{1}{|c|}{} \\
\cline{1-3}
&&&\\
\cline{1-4}
\end{tabular}\\
\noalign{\vskip2mm}
\hskip4mm \textrm{fig. 1}\\
\end{array}$ \hskip3mm
$
\begin{array}{ll}
\begin{tabular}{|cccc|l}
\cline{1-4}
\multicolumn{2}{|c|}{}&$-$ &+  &$\ \ell \ 6$ \\
\cline{1-2}
+ &+&$-$&$-$ & $\ \ell \ 7$\\
$-$ &+&+&$-$  &\ $\ell\ 8$     \\
+&+&+&+  &$\ \ell \ 9$ \\
+&$-$&+&$-$   &\ $\ell \ 10$\\
\cline{3-4}
$-$ &$-$ &\multicolumn{2}{|c|}{}   &\ $\ell\ {11}$\\
\cline{1-4}
\end{tabular}\\
 &\hskip-1cm S(42)-S(22) \\
\hskip4mm \textrm{fig. 2} &\\
 &\\
  &\\
   &\\
    &\\  &\\  &\\  &\\  &\\ &\\
\end{array}
$
\end{center}


\vskip2mm

2) $16 \leq n-m\leq 64$. It suffices now to consider the first 256 terms (the first 64 lines of the table). The idea in order to pick the infimum is to start from $S(4\times 11)-S(4\times 5)$ and the analogues, and modify the figure in order to diminish $||S(n)-S(m)||^2$ (fig. 2). As a first example, $S(44)-S(20)=2u_2+2u_3$ (obtained by replacing $u_1+u_2+u_3+u_4$ and $u_4$ by $u_1-u_2-u_3+u_4$ in the expression of $S(11) -S(5)$), and the modification provides $S(42)-S(22)= u_1+u_2+u_3-u_4$. The result is
$$
\inf ||S(n)-S(m)||^2=4
$$
realized for $(22,42)$ and $(214,234)$, with 
$$
\begin{array}{cc}
  S(42)-S(22) &=u_1 +u_2-u_3-u_4\,,\\
  S(234)-S(214)&=-u_1-u_2-u_3+u_4\,.   \\  
\end{array}
$$


This proves Lemma 2 with $A=16$.

\vskip2mm

Actually the proof can be given in a more concentrated form. It is enough to show that $||S(n)-S(m)||^2 \leq 3$ is impossible when $n-m\geq 16$. Let us assume \textit{ab absurdo} that $||S(n) -S(m)||^2\leq 3$.  Let us add or remove the minimal number of terms in order to transform $S(n)-S(m)$ into a difference of the form $S(4n') - S(4m')$ (that is, to transform the figure $S(n)-S(m)$ into a rectangle). In general, this minimal number is $\leq 4$ and has the same parity as $||(S(n)-S(m)||^2$ ; here it is $\leq 3$ and the resulting $S(4n')-S(4m')$ has its squared norm $\leq 2$, therefore $||S(n') -S(m')||^2\leq 3$ and the process goes on until we reach $S(n^{\times}) - S(m^{\times})$ with $n^{\times}\leq 64$. Then we know the possible pairs $(m^{\times},n^{\times})$, namely (5,11), (23,29), (35,41) and (53,59), and the reverse process never gives a squared norm $\leq 3$.

\vskip2mm

\noindent \textbf{5.6} Remarks and questions

The estimates we gave for $b$ and $a$ are quite rough. We can ask for better estimates and conjectures. The actual problem, of a combinatorial or arithmetical nature, is to compute these numbers exactly.

We were interested in estimating $b$ from above and $a$ from below. Examples provide estimates in the opposite direction :
$$
\begin{array}{cc}
  b &  \geq \frac{||S(17)||}{\sqrt{17}} = \frac{5}{\sqrt{17}}  \geq 1.21\hfill  \\
  \\
 a & \leq \frac{||S(42)-S(22)||}{\sqrt{20}} = \frac{2}{\sqrt{20}} = \frac{1}{\sqrt{5}} \leq 0.45   \\  
\end{array}
$$

It seems not impossible that the estimate for $\alpha$ is precise, that is $a=\frac{1}{\sqrt{5}}$. A careful investigation of the table would confirm or disprove this conjecture. It would lead also to a better estimate for $b$.

\section{Projections of the curve}

\noindent\textbf{6.1} \hskip2mm The direction of $u_0$ is special : all first coordinates of the $S(n) $ are $\geq 0$. That means that the partial sums $S_0(n)$ of the original series described in 3.0 are positive.

A simple way to see it is to use Lemma 1 (of 5.2). Since $S_0(n)\geq 1$ for $n=1,2,3,4,5,6,7,8$, we have $S_0(t)>0$ for $\frac{1}{2}\leq t \leq 8$, therefore (changing $t$
 into $16^k t$), $S_0(t)>0$ for all $t>0$.
 
 \vskip4mm

\noindent\textbf{6.2} \hskip2mm We are mainly interested in the three--dimensional projections of the curve. It seems likely that all parallel projections of the curve on a three--dimensional subspace of $\bbbr^4$ have an infinity of double points. The question can be formulated in the equivalent forms :

\vskip2mm

1) is every direction in $\bbbr^4$ the direction of some $S(t)-S(s)$ ?

2) are the $\frac{S(n)-S(m)}{\sqrt{n-m}}$ $(n>m)$ dense on the sphere $S^3$ ?

\vskip4mm

\noindent\textbf{6.3} \hskip2mm Let us project the curve from $0$ on the sphere $S^3$, that is consider 
$$
\cals(t)=\frac{S(t)}{||S(t)||} \qquad (t>0)\,.
$$
We obtain a closed $\calc$, image of any interval $[a,16a]$ by $\cals(\cdot)$.

$\calc$ is invariant under the isometries of $\bbbr^4$ defined by $\frac{1}{2}M$ and $\frac{1}{\sqrt{2}}T$ (see 4.1 and 4.2). The first takes the form
$$
\frac{1}{2}M' =
\left(
\begin{array}{cccc}
  1\ &0   &0 &0   \\
  0\ &  1 & 0 &0   \\
  0\ &  0 &-1 &0 \\
  0\ &0&0&-1 
\end{array}
\right)
$$
with respect to the orthonormal basis $(v_0,v_1,v_2,v_3)$ and the second
$$
\frac{1}{\sqrt{2}}T' =
\left(
\begin{array}{cccc}
  1&0   &0 &0   \\
  0&  -1 & 0 &0   \\
  0&  0 &0 &-1 \\
  0 &0&1&0 
\end{array}
\right)
$$
with respect to the orthonormal basis $(w_0,w_1,w_2,w_3)$. The vectors $v_0$ and $v_1$ ($w_0$ and $w_1$ as well) generate a plane $P$ such that the mapping $\cals(t)\lgr \cals(4t)$ is an orthogonal symmetry with respect to $P$. For the projection of $\calc$ on $P$, the change of $t$ into $2t$ means a symmetry with respect to the line generated by $w_0$.

$\calc$ has a double point at $t={1\over3}$, $t'={4\over3}$ : $\calt({4\over3})=\calt({1\over3})$.

In order to prove it, we expand $t$ and $t'$ in base $4$ (we underline the expansion)
$$
t= \underline{0. 1111}\cdots \qquad   t'= \underline{1. 1111}\cdots \,.
$$
Using base $4$ again, we easily obtain the figures and the values of $S(\underline{1})$, $S(\underline{11})$, $S(\underline{111})$ and so on :  $S(1) = u_0=(1\ 0\ 0\ 0)$
$$
\begin{array}{ll}
 S(\underline{11}) \hfill   &=S(\underline{10})+ (S(\underline{11}) - S(\underline{10}))= (1\ 1\ 1\  1)+(1\ 0\ 0\ 0) \hfill \\
 \noalign{\smallskip}
 S(\underline{111})    &=S(\underline{100})+ (S(\underline{110}) - S(\underline{100})) +S(\underline{111}) - ( S(\underline{110}) \hfill  \\ 
  \noalign{\smallskip}
 &=(4\ 0\ 0\ 0) + (1\ 1\ 1\ 1) - (1\ 0\ 0\ )\hfill\\
  \noalign{\smallskip}
 S(\underline{1111}) &= S(\underline{1000}) +(S(\underline{1100}) - S(\underline{1000})) + (S(\underline{1110}
) - S(\underline{1100}))\\
&\hskip6cm + (S(\underline{1111}) - S(\underline{1110}))\\
 \noalign{\smallskip}
&=(4\ 4\ 4\ 4) + (4\ 0\ 0\ 0) - (1\ 1\ 1\ 1) + (1\ 0\ 0\ 0)\hfill\\
 \noalign{\smallskip}
S(\underline{11111})  &= (16\ 0\ 0\ 0)+(4\ 4\ 4\ 4) - (4\ 0\ 0\ 0)+ (1\ 1\ 1\ 1) - (1\ 0\ 0\ 0)\\
 \noalign{\smallskip}
S(\underline{111111}) &= (16\ 16\ 16\ 16) + (16\ 0\ 0\ 0) - (4\ 4\ 4\ 4) + (4\ 0\ 0\ 0)\\
&\hskip6cm - (1\ 1\ 1\ 1) + (1\ 0\ 0\ 0)
\end{array}
$$

The ratio between two consecutive vectors tend to 2 (meaning that the ratios of coordinates tend to 2), hence
$$
\cals(t') =\calt(t) \qquad \big(t = \frac{1}{3}\big)\,.
$$

By isometry we also have
$$
\calt(2t') = \calt(2t)\,.
$$

These double points are contained in the plane $P$, and they are symmetric with respect to the line generated by~$w_0$.

I believe, but did not prove, that these are the only multiple points of the curve $\calc$. In that case, $\calc$ is a Brownian quasi--helix (actually, a Brownian quasi--circle) on some $4$--covering of the sphere~$S^3$.

\vskip4mm

\noindent\textbf{6.4} \hskip2mm One can see the curve $\calc$ in two other ways

First, taking into account that the first coordinate $S_0(t)$ is always positive, we can consider
$$
\calr(t) = \frac{S(t)}{S_0(t)}
$$
and the curve $\calc'$ described by $\calr(\cdot)$, projection of the original curve with a source at $0$ and a screen at the hyperplane $x_0=1$.

Symmetries and double points can be studied on this model as well.

Secondly, we obtain a projective model of $\calc$, say, $\calc''$, on choosing four points $A_0,A_1,A_2,A_3$ in $\bbbr^4$, defining $A_{j+4}=A_j$ $(j=0,1,\ldots)$, starting with a point $M_0=A_0$ and defining the sequence of points
$$
M_{n+1} = \frac{1}{a_0+ a_1+\cdots a_n} ((a_0+a_1+\cdots a_{n-1}) M_n +a_n A_n)\,.
$$
 
 \vskip2mm
 
 Some real figures would help. If a reader is willing  to draw figures of the above curves, I'll appreciate to see them.

\vskip0.8cm

\hfill \textsf{jean-pierre.kahane@math.u-psud.fr}

\end{document}